\documentclass[11pt]{article}
\usepackage[margin=1in]{geometry}
\usepackage{amsmath,amssymb,amsthm}
\usepackage{tikz}
\usetikzlibrary{calc}
\usepackage{hyperref}
\hypersetup{colorlinks=true,linkcolor=blue,citecolor=blue,urlcolor=blue}

\newtheorem{theorem}{Theorem}

\newtheorem{corollary}{Corollary}
\theoremstyle{remark}

\newcommand{\Aff}{\operatorname{Aff}}
\newcommand{\SL}{\operatorname{SL}}

\title{Veech dichotomy for specific covering surface over Veech surface}
\author{Haoran Hu}
\date{\today}

\begin{document}
\maketitle

\begin{abstract}
It is known that Veech dichotomy, which states that minimal directions are uniquely ergodic, holds for Veech surface and covering surface branched at a single non-periodic connection point of Veech surface. We give a sufficient condition for Veech dichotomy holds for a finite translation cover of a Veech surface, and applies it on branched cover of square-tiled surfaces to construct examples that satisfy Veech dichotomy but aren't Veech surfaces.
\end{abstract}

\section{Introduction}
In 1989, Veech proved that a translation surface whose affine derivative group is a lattice in $\SL(2,\mathbb R)$ satisfies a strong dynamical dichotomy: For every direction, the direction flow is either completely periodic or minimal and uniquely ergodic\cite{Veech}, which is now also called Veech dichotomy. This raised the question of whether the Veech dichotomy characterizes lattice surfaces.

Branched coverings provided a way to separate these geometric and dynamical properties. In 2004, Hubert and Schmidt constructed covers of lattice surfaces with infinitely generated Veech groups \cite{HS}. Building on their construction, in 2008, Smillie and Weiss proved that the converse to Veech's theorem fails: there are non-Veech translation surfaces satisfying the Veech dichotomy \cite{SW}. Their examples use covers branched over a single nonperiodic connection point. More generally, their compactness argument, together with Masur's recurrence criterion, gives unique ergodicity in every minimal direction for a finite cover of a Veech surface branched over one value.

The role of the branch locus is further illustrated by Cheung, Hubert and Masur \cite{CHM}. In 2008 they constructed branched covers satisfying \emph{topological dichotomy}---every direction is completely periodic or minimal---but having minimal directions that are not uniquely ergodic. Thus topological dichotomy and \emph{strict ergodicity}, meaning unique ergodicity in every minimal direction, are distinct requirements. The author studies sufficient conditions for strict ergodicity when there is more than one branch value.

On a Veech surface, a minimal direction gives a Teichm\"uller ray with infinitely many returns to a compact set. For a branched cover, the covering surface depends also on the positions of the branch values. Now we ask when the same times give compact returns for the cover. Our condition concerns the movement of these values under the affine automorphisms of the base.

\section{A criterion for branched covers}
Write $\Aff^+(X,\omega)$ for the orientation-preserving affine automorphisms of a compact translation surface. The surface is Veech if the group of their derivatives is a lattice in $\SL(2,\mathbb R)$. A directional flow is minimal if every complete nonsingular trajectory is dense. Unique ergodicity means that normalized flat area is the unique invariant probability measure for the nonsingular flow. For an ordered tuple of $m$ points, put
\[
 \Delta_m=\bigcup_{i<j}\{x_i=x_j\}\subset X^m,\qquad
 \operatorname{Conf}_m(X)=X^m\setminus\Delta_m .
\]
The affine group acts on an ordered tuple by applying the same automorphism to each coordinate. For an unbranched cover we take $\operatorname{Conf}_0(X)$ to consist of the empty tuple and $\Delta_0=\varnothing$.

\begin{theorem}[Compact branch configurations]\label{thm:configuration}
Let $(X,\omega)$ be a compact Veech translation surface of genus $g\geq2$. Let $\pi:Y\to X$ be a connected finite-degree translation cover with $\eta=\pi^*\omega$, and let $\mathbf b=(b_1,\ldots,b_m)$ be its ordered tuple of distinct branch values. Suppose there is a compact $\Aff^+(X,\omega)$-invariant set
\[
 F\subset\operatorname{Conf}_m(X)
\]
containing $\mathbf b$. Then every minimal direction on $Y$ is uniquely ergodic.
\end{theorem}

Equivalently, it is sufficient that the closure of $\Aff^+(X,\omega)\cdot\mathbf b$ in $X^m$ miss the fat diagonal. Such a closure is a compact invariant subset of $\operatorname{Conf}_m(X)$. For $m\le1$, the condition is automatic.

\subsection{The space of covers and its projections}
Let $\mathcal M_g$ be the moduli space of smooth compact Riemann surfaces of genus $g$. We use the following conventions for differentials and markings.
\begin{description}
\item[$\mathcal Q_g$.] A point is a pair $(C,q)$ consisting of a curve $C$ in $\mathcal M_g$ and a holomorphic quadratic differential $q\in H^0(C,K_C^{\otimes 2})$ with $\int_C|q|=1$. Two pairs represent the same point when a biholomorphism pulls the second differential back to the first. The topology is that of the unit-area bundle over $\mathcal M_g$. Zeros may merge as the differential varies on a smooth curve.
\item[$\mathcal Q_{g,b}$.] A point is a pair $(C,q)$ as above together with $b$ ordered, pairwise distinct points $\mathbf x=(x_1,\ldots,x_b)$ on $C$. Thus $\mathcal Q_{g,0}=\mathcal Q_g$. The biholomorphisms in the equivalence relation preserve the labels. In this paper we take $q$ holomorphic on the whole curve. A mark may lie at a zero of $q$ \cite[p.~1966]{SW}. The map $\bar\Phi_{g,b}$ forgetting $\mathbf x$ has fiber $\operatorname{Conf}_b(C)$ over a fixed representative, modulo the finite automorphism group of $(C,q)$ on coarse moduli.
\item[$\mathcal Q(\pi_0)$.] Fix a reference topological branched cover $\pi_0:\widetilde S\to S$ of degree $d$, with base genus $g$, covering genus $h$, and $b$ labeled branch values. A point of $\mathcal Q(\pi_0)$ records $(C,q,\mathbf x)$ together with a holomorphic branched cover $p:D\to C$ of that topological type. This means that there are orientation-preserving diffeomorphisms $j:S\to C$ and $\tilde j:\widetilde S\to D$ respecting the branch labels and satisfying $p\circ\tilde j=j\circ\pi_0$. We identify two such records by pairs of biholomorphisms commuting with the covering maps and preserving the marked base differential. The covering differential is determined by these data. We give it area one by taking $q_D=d^{-1}p^*q$, since $\int_D|p^*q|=d$.
\end{description}
For a given translation cover $\pi$, the notation $\mathcal Q(\pi)$ refers to this space with $\pi_0$ the underlying topological map of $\pi$. Writing $u=[p:D\to C,q,\mathbf x]$, we have the maps
\[
 \mathcal Q_{g,b}\ \xleftarrow{\ P\ }\ \mathcal Q(\pi_0)\
 \xrightarrow{\ R\ }\ \mathcal Q_h\
 \xrightarrow{\ \operatorname{pr}_h\ }\ \mathcal M_h .
\]
In this diagram $P(u)=[C,q,\mathbf x]$, $R(u)=[D,q_D]$, and $\operatorname{pr}_h$ retains the underlying curve. Simultaneously deforming the flat charts on base and cover gives the $\SL(2,\mathbb R)$ action on $\mathcal Q(\pi_0)$ and makes $P$ and $R$ equivariant.

\paragraph{Properness over the marked base.}
Over a fixed pointed base, a connected degree-$d$ cover is described by a transitive monodromy representation of $\pi_1(C\setminus\{x_1,\ldots,x_b\})$ into $S_d$, up to relabeling the sheets. The cycles around the punctures specify local ramification. There are finitely many choices since the group is finitely generated and $S_d$ is finite. As the distinct marks and the base curve vary in a small moduli chart, local identifications of the punctured curves transport these choices. This gives the local sheets of the Hurwitz space for the chosen reference type. Smillie and Weiss establish continuity and equivariance in Proposition~6 and properness of $P$ in Proposition~7. Corollary~17 identifies $P$ as a finite orbifold covering \cite[pp.~1962, 1970]{SW}. In particular, $P^{-1}(B)$ is compact whenever $B\subset\mathcal Q_{g,b}$ is compact.

For comparison, $\bar\Phi_{g,b}:\mathcal Q_{g,b}\to\mathcal Q_g$ forgets the branch marks themselves. Proposition~8 of \cite{SW} gives its properness for $b=1$ and $g\ge2$. For $b\ge2$, two marks on a fixed curve can converge to the same point, escaping compact subsets of the fiber while their unmarked base remains constant \cite[p.~1963]{SW}. The hypothesis of Theorem~\ref{thm:configuration} will give a compact set of marked bases to which the properness of $P$ applies.

\subsection{Recurrence and the proof}

We recall the base recurrence used in the proof. Let $G=\SL(2,\mathbb R)$ and let $\Gamma\subset G$ be the derivative group of a Veech surface. Given a direction $\theta$, choose a rotation $r_\theta\in G$ taking it to the vertical direction and set
\[
 g_t=\begin{pmatrix}e^t&0\\0&e^{-t}\end{pmatrix},\qquad H_t=g_t r_\theta.
\]
The ray $H_tX$ contracts the tested direction. If this direction is minimal, Veech recurrence gives a compact subset $C\subset G/\Gamma$ visited by $H_t\Gamma$ for arbitrarily large $t$ \cite[Proposition~5]{SW}.

Let $\rho:G\to G/\Gamma$ be the quotient map. It is open because $\rho^{-1}(\rho(U))=\bigcup_{\gamma\in\Gamma}U\gamma$ is open for every open $U\subset G$. For each $x\in C$, choose a lift $g_x$ and a relatively compact open neighborhood $U_x\subset G$. The sets $\rho(U_x)$ cover $C$. Choose a finite subcover $\rho(U_{x_1}),\ldots,\rho(U_{x_k})$. Then
\[
 K=\bigcup_{j=1}^k\overline{U_{x_j}}
\]
is compact and $C\subset\rho(K)$. For return times $t_n\to\infty$, choose $A_n\in K$ with $\rho(A_n)=\rho(H_{t_n})$. Since the quotient is by right cosets, $A_n=H_{t_n}\gamma_n$ for some $\gamma_n\in\Gamma$. 

\begin{proof}[Proof of Theorem~\ref{thm:configuration}]
Fix a direction $\theta$ that is minimal on $Y$. Its projected direction on $X$ cannot be periodic. In a periodic direction on the Veech surface $X$, the periodic decomposition has finitely many cylinders and boundary saddle connections. Each branch value lies on at most one periodic leaf, so a periodic leaf can be chosen to avoid the finite branch locus. The return around this leaf permutes finitely many sheets, so each lift closes after finitely many circuits. This contradicts minimality on $Y$. The projected direction is therefore nonperiodic, and Veech dichotomy implies it is minimal on $X$.

Apply the recurrence above to the base direction. We obtain $t_n\to\infty$ and $\gamma_n\in\Gamma$ such that $A_n=H_{t_n}\gamma_n\in K$. Choose affine automorphisms $\phi_n$ of $X$ with derivative $\gamma_n$. The branch tuple normalized along the return is
\[
 \mathbf b_n=\phi_n^{-1}(\mathbf b).
\]
By affine invariance, $\mathbf b_n\in F\subset X^m$. To keep track of the marks under deformation, let $d_A:X\to AX$ denote the natural deformation map for $A\in G$. This map identifies the underlying topological surfaces and has derivative $A$ in flat charts. We apply it coordinatewise to tuples. In particular, a tuple $\mathbf x$ on $X$ becomes $d_A(\mathbf x)$ on $AX$; writing $(AX,\mathbf x)$ is shorthand for this identification.

The actual branch tuple on the base at time $t_n$ is $d_{H_{t_n}}(\mathbf b)$. Consider the map
\[
 \widehat\phi_n
 =d_{H_{t_n}}\circ\phi_n\circ d_{A_n}^{-1}
 :A_nX\longrightarrow H_{t_n}X.
\]
Its derivative in flat charts is
\[
 D\widehat\phi_n=H_{t_n}\gamma_n A_n^{-1}=I,
\]
so it is a translation equivalence. It also sends the marks to the correct points:
\[
 \widehat\phi_n\bigl(d_{A_n}(\mathbf b_n)\bigr)
 =d_{H_{t_n}}\bigl(\phi_n(\mathbf b_n)\bigr)
 =d_{H_{t_n}}(\mathbf b).
\]
Thus the marked base at time $t_n$ is represented by
$(A_nX,q_{A_nX},d_{A_n}(\mathbf b_n))$, where $q_{AX}$ denotes the normalized square of the deformed form.

Here $F$ controls the tuples on the fixed surface $X$: the assertion is $\mathbf b_n\in F$, and we do not assert that the deformed tuple $d_{A_n}(\mathbf b_n)$ belongs to $F$. To obtain a compact set of marked surfaces, define
\[
 \Phi:K\times F\longrightarrow\mathcal Q_{g,m},\qquad
 \Phi(A,\mathbf x)=[AX,q_{AX},d_A(\mathbf x)],
 \qquad B:=\Phi(K\times F).
\]
The map $\Phi$ is continuous, so $B$ is compact as the image of the compact product $K\times F$. Each $d_A$ is a homeomorphism, so distinct marks remain distinct and $B$ lies in the distinct-mark locus. Since $(A_n,\mathbf b_n)\in K\times F$, every normalized marked base belongs to $B$.

Take $\pi_0$ to be the underlying topological map of $\pi$. Let $u_0\in\mathcal Q(\pi_0)$ represent the original cover with the unit-area base differential $q_X=\omega^2/\int_X|\omega^2|$ and covering differential $q_Y=d^{-1}\pi^*q_X$, where $d=\deg\pi$. Put $h=g(Y)$ and $u_n=H_{t_n}u_0$. The marked bases satisfy $P(u_n)=\Phi(A_n,\mathbf b_n)\in B$. Since $P$ is proper, $L=P^{-1}(B)$ is compact. Thus $u_n\in L$. The map $R$ is continuous and equivariant under the $G$-action, so
\[
 H_{t_n}(Y,q_Y)=R(u_n)\in R(L),
\]
where $R(L)$ is compact in $\mathcal Q_h$. It follows that
\[
 [H_{t_n}Y]\in C_h:=\operatorname{pr}_h(R(L))\subset\mathcal M_h,
\]
and $C_h$ is compact. This proves recurrence of the Teichm\"uller ray of $Y$ in unmarked genus-$h$ moduli. The argument uses the affine automorphisms of $X$ to choose representatives of $P(u_n)$. Properness then controls the covering data over all these representatives, whether or not a given base automorphism lifts to a self-map of $Y$.

The direction on $Y$ was assumed minimal. Masur's criterion says that a minimal directional flow whose Teichm\"uller ray is recurrent in $\mathcal M_h$ is uniquely ergodic \cite{Masur92}. Hence the flow in direction $\theta$ is uniquely ergodic. \end{proof}

\section{A torsion criterion on the period torus}
Suppose the absolute periods $\int_\gamma\omega$, for $\gamma\in H_1(X;\mathbb Z)$, form a lattice $\Lambda\subset\mathbb R^2$. Integrating $\omega$ from a base point defines a map $f:X\to E=\mathbb R^2/\Lambda$. A change of base point translates $f$ by a constant, so the difference $f(y)-f(x)$ is intrinsic. Write $E[N]=\{t\in E: Nt=0\}$ for the finite subgroup of $N$-torsion points.

\begin{corollary}[Torsion-separated branch values]\label{cor:torsion}
Let $(X,\omega)$ be a Veech translation surface of genus at least two whose
full absolute period group is a rank-two lattice $\Lambda$. Let
$f:X\to E=\mathbb R^2/\Lambda$ be the canonical period-torus map obtained
by integration of $\omega$, and fix $N\geq2$. Suppose that an ordered tuple
of distinct branch values $\mathbf b=(b_1,\ldots,b_m)$ satisfies
\[
 f(b_j)-f(b_i)\in E[N]\setminus\{0\}\qquad(i\ne j).
\]
Then every connected finite-degree translation cover of $X$ with precisely
these branch values has uniquely ergodic flow in
every minimal direction. If at least one branch value has infinite affine
orbit, the cover is not a Veech surface \cite[Proposition~2.1(c)]{CW}.
\end{corollary}
\begin{proof}
For an affine automorphism $\phi$ with derivative $A$, preservation of the
full absolute period group gives $A\Lambda=\Lambda$. Integration yields
$f\circ\phi=A f+c$ for a constant $c\in E$, so the constant cancels in
differences. The torsion conditions are therefore affine-invariant. The
set of tuples satisfying them is closed in $X^m$, hence compact. Since all
the required differences are nonzero, it avoids the fat diagonal. It is an
invariant compact subset of $\operatorname{Conf}_m(X)$ containing
$\mathbf b$, so Theorem~\ref{thm:configuration} applies. As $N$ varies, the two-point case covers every nonzero torsion difference.
\end{proof}

\section{Genus-four double covers of a three-square surface}
Take three unit squares labeled $1,2,3$. Glue the right edge of square $i$
to the left edge of square $h(i)$, and the top edge of square $i$ to the
bottom edge of square $v(i)$, by translations, where
\[
 h=(12),\qquad v=(13).
\]
The permutations act transitively, so the resulting surface $X$ is
connected. Their commutator is a $3$-cycle, giving a single vertex class
with cone angle $6\pi$. Thus $X$ has one zero of order two and genus two.
The square tiling defines a degree-three map
\[
 f:X\longrightarrow\mathbb R^2/\mathbb Z^2,
 \qquad f^*(dz)=\omega,
\]
branched only over the origin. Every absolute period is integral, while the
horizontal loop through square $3$ and vertical loop through square $2$
have periods $(1,0)$ and $(0,1)$. The full absolute period lattice is
therefore $\mathbb Z^2$. The derivative group is contained in
$\operatorname{SL}(2,\mathbb Z)$, and has finite index there because the
$\operatorname{SL}(2,\mathbb Z)$ orbit of a degree-three square-tiled
surface is finite. Hence $X$ is a Veech surface.

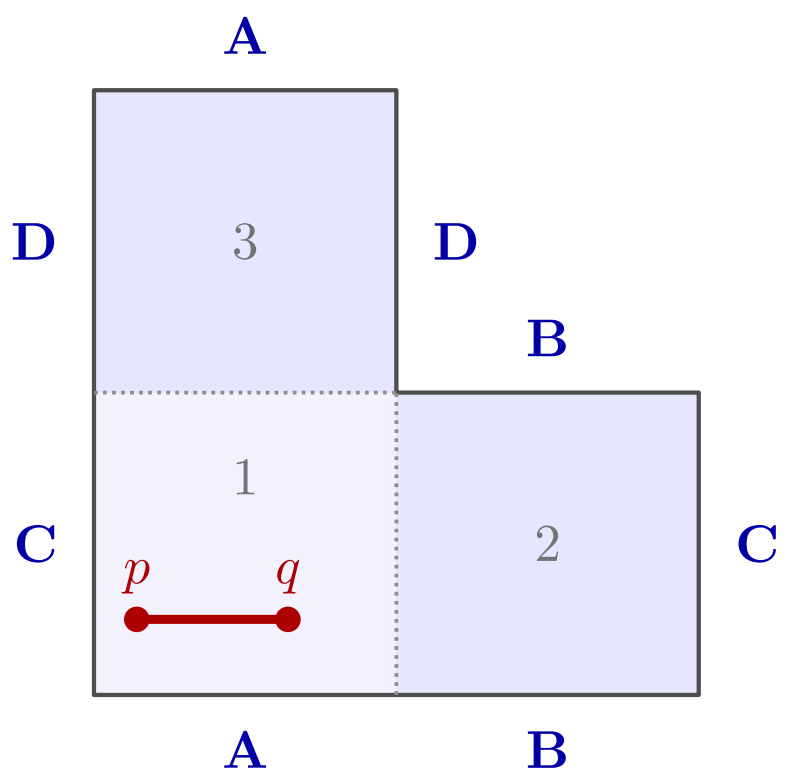
\begin{figure}[!htb]
\centering
\begin{tikzpicture}[x=2.0cm,y=2.0cm,line join=round,
  tile/.style={draw=black!70,line width=0.8pt},
  seam/.style={draw=black!45,densely dotted,line width=0.7pt},
  edge label/.style={font=\small\bfseries,text=blue!65!black},
  tile number/.style={font=\small,text=black!55},
  slit/.style={draw=red!68!black,line width=1.7pt},
  origami branch/.style={circle,fill=red!68!black,inner sep=1.7pt}]
  \fill[blue!5] (0,0) rectangle (1,1);
  \fill[blue!10] (1,0) rectangle (2,1);
  \fill[blue!10] (0,1) rectangle (1,2);
  \draw[tile] (0,0)--(2,0)--(2,1)--(1,1)--(1,2)--(0,2)--cycle;
  \draw[seam] (1,0)--(1,1);
  \draw[seam] (0,1)--(1,1);
  \node[tile number] at (.5,.72) {1};
  \node[tile number] at (1.5,.5) {2};
  \node[tile number] at (.5,1.5) {3};
  \node[edge label,below=3pt] at (.5,0) {A};
  \node[edge label,above=3pt] at (.5,2) {A};
  \node[edge label,below=3pt] at (1.5,0) {B};
  \node[edge label,above=3pt] at (1.5,1) {B};
  \node[edge label,left=3pt] at (0,.5) {C};
  \node[edge label,right=3pt] at (2,.5) {C};
  \node[edge label,left=3pt] at (0,1.5) {D};
  \node[edge label,right=3pt] at (1,1.5) {D};
  \coordinate (p) at ({sqrt(2)/10},.25);
  \coordinate (q) at ({sqrt(2)/10+.5},.25);
  \draw[slit] (p)--(q);
  \node[origami branch] at (p) {};
  \node[origami branch] at (q) {};
  \node[font=\small\bfseries,text=red!68!black,above=1pt] at (p) {$p$};
  \node[font=\small\bfseries,text=red!68!black,above=1pt] at (q) {$q$};
\end{tikzpicture}
\caption{The three-square base. Boundary segments with the same letter are glued by translation. The points $p=(\sqrt{2}/10,1/4)$ and $q=p+(1/2,0)$ illustrate one member of Corollary~\ref{thm:family}. Cutting along the marked segment and cross-gluing two copies gives one of its double covers.}
\label{fig:three-square}
\end{figure}

\paragraph{The two-parameter family.}
In square $1$, choose
\[
 0<a<\tfrac12,\qquad 0<b<1,\qquad
 p_{a,b}=(a,b),\quad q_{a,b}=(a+\tfrac12,b).
\]
The points are regular and distinct, and their period-torus difference is
$(1/2,0)$, a nonzero two-torsion point. Figure~\ref{fig:three-square} shows a representative slit in square~1.

The set of connected double
covers simply branched exactly over these labeled points has 16 elements up to isomorphism over the identity of the marked base.
Indeed, the monodromy character
\[
 \chi:\pi_1(X\setminus\{p_{a,b},q_{a,b}\})\to\mathbb Z/2\mathbb Z
\]
must take value $1$ on each puncture loop. In the standard presentation
\[
 \langle a_1,b_1,a_2,b_2,c_p,c_q\mid
 [a_1,b_1][a_2,b_2]c_pc_q=1\rangle,
\]
the four handle values are free, giving $2^4=16$ characters. Each cover is
connected because the puncture monodromy is nontrivial. Pulling back the complex structure gives a holomorphic cover of the punctured curve. At each puncture, filling in the unique ramification point with the local model $w\mapsto w^2$ gives its unique extension up to isomorphism over the base. The slit joining
$p_{a,b}$ to $q_{a,b}$ gives one of these covers. Its monodromy is mod-two
intersection with the slit, so handle loops chosen disjoint from the slit
have value zero and the two puncture loops have value one.

\begin{corollary}[The two-parameter family]\label{thm:family}
For every $(a,b)$ in the indicated open rectangle, every connected double
cover $Y_{a,b}\to X$ branched exactly at $p_{a,b},q_{a,b}$ has genus four, zero
orders $(2,2,1,1)$, and uniquely ergodic flow in every minimal direction.
Moreover, $Y_{a,b}$ is a Veech surface if and only if $a,b\in\mathbb Q$.
\end{corollary}
\begin{proof}
Riemann--Hurwitz gives
\[
 2g(Y_{a,b})-2=2(2g(X)-2)+2=6,
\]
so $g(Y_{a,b})=4$. The unique order-two zero on $X$ has two unbranched
preimages of order two. Each of the two simple branch points is regular on
$X$ and contributes a simple zero of the pulled-back differential. Hence
the zero orders are $(2,2,1,1)$. Corollary~\ref{cor:torsion}
gives unique ergodicity in every minimal direction for every such cover.

If $a,b$ are rational, subdivide the three base squares by a sufficiently
fine rational grid for which both branch points are vertices. Its pullback
to any such cover is a square tiling after a common rescaling. Thus
$Y_{a,b}$ is square-tiled and Veech.

For the converse, consider the segments in square $1$ from its lower corners
$(0,0)$ and $(1,0)$ to $p_{a,b}$. The corners represent the unique zero of
$X$. The segment interiors lie below height $b$ and avoid the branch points.
Their lifts are saddle connections with holonomies $(a,b)$ and $(a-1,b)$,
and slopes
\[
 r=\frac{b}{a},\qquad s=\frac{b}{a-1}.
\]
If both slopes were rational, then $a=s/(s-r)$ and $b=ra$ would both be
rational. Thus for $(a,b)\notin\mathbb Q^2$, at least one of these saddle
connections has irrational slope. A closed regular geodesic on the cover
projects under $f$ composed with the covering map to a closed straight
geodesic on the square torus, so it has rational slope. The irrational
saddle-connection direction contains no cylinder, contradicting complete
periodicity of saddle-connection directions on a Veech surface
\cite{Veech}. Hence the cover is not Veech.

\paragraph{The branch points at compact returns.}
In this family the separation used in Theorem~\ref{thm:configuration} can be read off from the square torus $T=\mathbb R^2/\mathbb Z^2$. For the normalized marks $p_n=\phi_n^{-1}(p_{a,b})$ and $q_n=\phi_n^{-1}(q_{a,b})$, the displacement $f(q_n)-f(p_n)$ is one of the three nonzero elements of $T[2]$. Their distances from the origin in the unit flat torus are all at least $1/2$. The map $f$ is length-nonincreasing, so
\[
 d_X(p_n,q_n)\ge \tfrac12,\qquad
 d_{A_nX}(p_n,q_n)\ge \tfrac12\min_{A\in K}\sigma_{\min}(A)>0 ,
\]
where $\sigma_{\min}$ is the least singular value and $K$ is the compact set of return representatives in the proof of Theorem~\ref{thm:configuration}.

An individual branch value can still approach the old order-two zero. In a local coordinate $u$ with $\omega=u^2\,du$, simple branching at $u=t$ has the model
\[
 u=w^2+t,\qquad \pi^*\omega=2w(w^2+t)^2\,dw .
\]
For $t\ne0$ the ramification point gives a simple zero at $w=0$, and the old zero has two order-two preimages at $w=\pm\sqrt{-t}$. At $t=0$ the limiting form is $2w^5\,dw$ on the smooth $w$-disk. The separation estimate keeps the second branch value away from this same base zero. Such a limit has orders $(5,1)$ on a smooth genus-four surface. The compactness statement used in Theorem~\ref{thm:configuration} concerns the underlying curve in $\mathcal M_4$ and accommodates this smooth limit.
\end{proof}
\section{Double covers of the regular double pentagon}
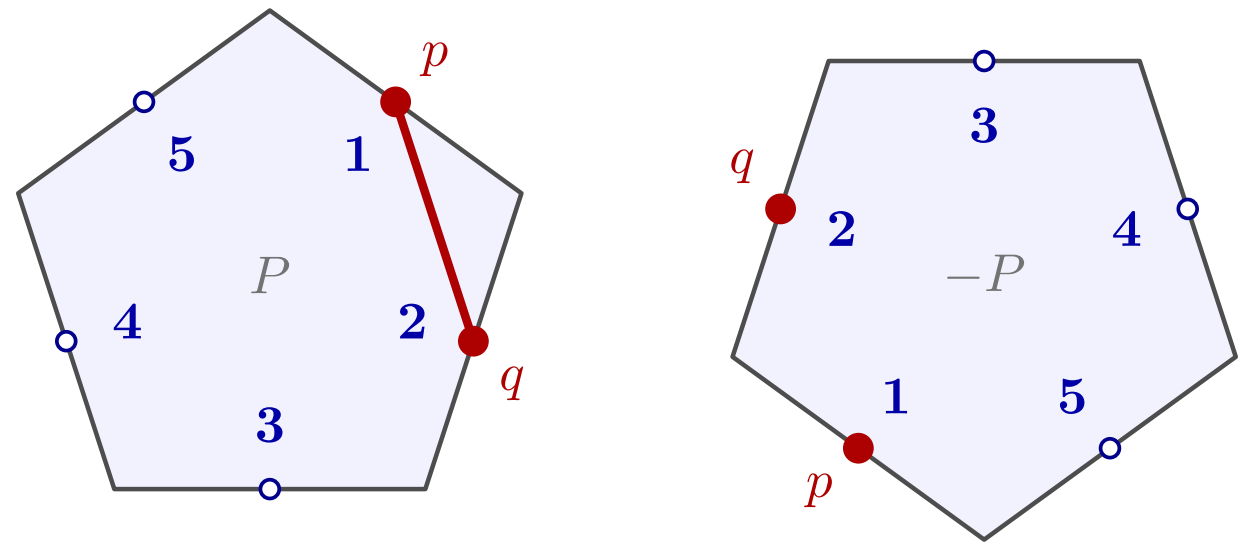
\begin{figure}[!htb]
\centering
\begin{tikzpicture}[x=1.75cm,y=1.75cm,line join=round,
  polygon rim/.style={draw=black!70,line width=0.85pt},
  edge label/.style={font=\small\bfseries,text=blue!65!black},
  paired mark/.style={circle,draw=blue!55!black,fill=white,line width=.7pt,inner sep=1.25pt},
  branch mark/.style={circle,draw=red!68!black,fill=red!68!black,line width=.7pt,inner sep=1.8pt},
  branch label/.style={font=\small\bfseries,text=red!68!black}]
  \coordinate (A0) at (-1.35,1);
  \coordinate (A1) at ({-1.35+cos(18)},{sin(18)});
  \coordinate (A2) at ({-1.35+cos(-54)},{sin(-54)});
  \coordinate (A3) at ({-1.35+cos(-126)},{sin(-126)});
  \coordinate (A4) at ({-1.35+cos(162)},{sin(162)});
  \coordinate (B0) at (1.35,-1);
  \coordinate (B1) at ({1.35-cos(18)},{-sin(18)});
  \coordinate (B2) at ({1.35-cos(-54)},{-sin(-54)});
  \coordinate (B3) at ({1.35-cos(-126)},{-sin(-126)});
  \coordinate (B4) at ({1.35-cos(162)},{-sin(162)});
  \coordinate (Cleft) at (-1.35,0);
  \coordinate (Cright) at (1.35,0);
  \fill[blue!5] (A0)--(A1)--(A2)--(A3)--(A4)--cycle;
  \fill[blue!5] (B0)--(B1)--(B2)--(B3)--(B4)--cycle;
  \draw[polygon rim] (A0)--(A1)--(A2)--(A3)--(A4)--cycle;
  \draw[polygon rim] (B0)--(B1)--(B2)--(B3)--(B4)--cycle;
  \node[font=\small,text=black!55] at (-1.35,0) {$P$};
  \node[font=\small,text=black!55] at (1.35,0) {$-P$};
  \path
    coordinate (mA1) at ($(A0)!.5!(A1)$)
    coordinate (mA2) at ($(A1)!.5!(A2)$)
    coordinate (mA3) at ($(A2)!.5!(A3)$)
    coordinate (mA4) at ($(A3)!.5!(A4)$)
    coordinate (mA5) at ($(A4)!.5!(A0)$)
    coordinate (mB1) at ($(B0)!.5!(B1)$)
    coordinate (mB2) at ($(B1)!.5!(B2)$)
    coordinate (mB3) at ($(B2)!.5!(B3)$)
    coordinate (mB4) at ($(B3)!.5!(B4)$)
    coordinate (mB5) at ($(B4)!.5!(B0)$);
  \foreach \i/\m in {1/mA1,2/mA2,3/mA3,4/mA4,5/mA5}
    \node[edge label] at ($ (\m)!.30!(Cleft) $) {\i};
  \foreach \i/\m in {1/mB1,2/mB2,3/mB3,4/mB4,5/mB5}
    \node[edge label] at ($ (\m)!.30!(Cright) $) {\i};
  \foreach \m in {mA3,mA4,mA5,mB3,mB4,mB5}
    \node[paired mark] at (\m) {};
  \foreach \m in {mA1,mB1,mA2,mB2}
    \node[branch mark] at (\m) {};
  \node[branch label,above right=1pt] at (mA1) {$p$};
  \node[branch label,below right=1pt] at (mA2) {$q$};
  \node[branch label,below left=1pt] at (mB1) {$p$};
  \node[branch label,above left=1pt] at (mB2) {$q$};
  \draw[draw=red!68!black,line width=1.6pt] (mA1)--(mA2);
\end{tikzpicture}
\caption{The double-pentagon base. Equal numbers pair parallel sides by translation. The two copies of each filled midpoint represent a single branch value, $p$ or $q$. The hollow midpoints represent the other three regular periodic points. The slit is drawn in $P$.}
\label{fig:pentagon}
\end{figure}

Let $X_D$ be the surface obtained by gluing each side of a regular pentagon $P$ to the parallel side of a second pentagon $-P$ by translation, as in Figure~\ref{fig:pentagon}. All ten vertices give one point $z$ of angle $6\pi$, so $X_D$ has genus two and one zero of order two. The side midpoints give five regular points $W=\{w_1,\ldots,w_5\}$ on the quotient. This surface is primitive and nonarithmetic, and its affine-periodic regular points are precisely the elements of $W$ \cite[Theorem~1.3, Lemma~2.7 and proof of Corollary~1.7]{ASZ}.

\begin{corollary}[The double pentagon]\label{cor:pentagon}
For distinct regular $p,q\in X_D$, the pair $(p,q)$ satisfies the compact-configuration hypothesis of Theorem~\ref{thm:configuration} if and only if $p,q\in W$. For any such pair, every connected double cover branched exactly at $p$ and $q$ lies in $\mathcal H(2,2,1,1)$, is a Veech surface, and has uniquely ergodic flow in every minimal direction.
\end{corollary}
\begin{proof}
For the pictured cover, take the two adjacent midpoint classes $p,q$ in Figure~\ref{fig:pentagon}, cut two copies of $X_D$ along the displayed segment, and glue opposite banks across the copies. This gives a connected double cover branched exactly at $p,q$.

The affine group permutes the five elements of $W$. For any distinct $p,q\in W$, the set $(W\times W)\setminus\Delta_2$ is therefore a finite invariant set of distinct pairs. Applying Theorem~\ref{thm:configuration} gives the claimed unique ergodicity. Periodicity of both branch values gives the lattice property by \cite[Proposition~2.1(c)]{CW}. Simple branching at the two regular points creates two simple zeros, and the double zero $z$ has two order-two preimages. Riemann--Hurwitz gives $2g(Y)-2=2(2\cdot2-2)+2=6$.

For necessity, suppose that at least one coordinate has infinite affine orbit. By Hubert, Schmoll and Troubetzkoy, the closure of the joint orbit in the ordinary topology of $X_D^2$ is a complex-algebraic set \cite[Proposition~3]{HST}. Since the orbit is infinite, this closure has a positive-dimensional component. Their classification applies because $X_D$ is primitive and nonarithmetic and its unique cone point has angle $2\pi\cdot3$, with $3$ prime \cite[Theorem~8]{HST}. A component of complex dimension two equals $X_D^2$. A component of complex dimension one is a coordinate fiber, the diagonal, or the graph of the hyperelliptic involution $\iota$. Each meets $\Delta_2$. A coordinate fiber through $r$ contains $(r,r)$, and the graph of $\iota$ contains $(w,w)$ at every Weierstrass point $w$. Hence an orbit closure avoiding $\Delta_2$ must be finite. Both coordinates are then regular periodic points and belong to $W$.
\end{proof}

\section*{Acknowledgments}

The author would like to sincerely thank Professor Yitwah Cheung for helpful suggestions and discussions on possibility to construct more general surfaces satisfying Veech dichotomy.

A Lean~4 formalization of
Theorem~\ref{thm:configuration}, using mathlib, is available
in the \href{https://github.com/Ryan-Hao-Jan/Veech-dichotomy-for-cover}
{accompanying code repository}, together with instructions
for reproducing the verification.

ChatGPT was used as an AI-assisted research tool in the development of this work. The tool was used to develop the proof, format \LaTeX{} manuscript and lean. The author have verified the proof independently and modified some necessary details. The author takes full responsibility for all statements and proofs in the paper.


\begin{thebibliography}{9}
\bibitem{ASZ} P.~Apisa, R.~M. Saavedra and C.~Zhang, Periodic points on the regular and double $n$-gon surfaces, \emph{Geometriae Dedicata} \textbf{216} (2022), article 69.

\bibitem{CW} M.~Cohen and B.~Weiss, Parking garages with optimal dynamics, \emph{Geometriae Dedicata} \textbf{161} (2012), 157--167.

\bibitem{CHM} Y.~Cheung, P.~Hubert and H.~Masur, Topological dichotomy and strict ergodicity for translation surfaces, \emph{Ergodic Theory Dynam. Systems} \textbf{28} (2008), no.~6, 1729--1748.

\bibitem{HS} P.~Hubert and T.~A. Schmidt, Infinitely generated Veech groups, \emph{Duke Math. J.} \textbf{123} (2004), no.~1, 49--69.

\bibitem{HST} P.~Hubert, M.~Schmoll and S.~Troubetzkoy, Modular fibers and illumination problems, \emph{Int. Math. Res. Not. IMRN} (2008), article rnn011.


\bibitem{Masur92}
H.~Masur,
Hausdorff dimension of the set of nonergodic foliations
of a quadratic differential,
\emph{Duke Math. J.} \textbf{66} (1992), no.~3, 387--442.

\bibitem{SW} J.~Smillie and B.~Weiss, Veech's dichotomy and the lattice property, \emph{Ergodic Theory Dynam. Systems} \textbf{28} (2008), no.~6, 1959--1972.

\bibitem{Veech} W.~A. Veech, Teichm\"uller curves in moduli space, Eisenstein series and an application to triangular billiards, \emph{Invent. Math.} \textbf{97} (1989), 553--584.
\end{thebibliography}
\end{document}